\documentclass[11pt]{article}
\usepackage{amsfonts}
\usepackage{color}
\newtheorem{theo}{Theorem}[section]
\newtheorem{lem}[theo]{Lemma}
\newtheorem{cor}[theo]{Corollary}

\newcommand{\mysection}[1]{\section{#1} \setcounter{equation}{0}}
\newcommand{\proof}{{\sc Proof.} \quad}
\newcommand{\proofc}{{\sc Proof} \ }
\newcommand{\be}{\begin{equation} \label}
\newcommand{\ee}{\end{equation}}
\newcommand{\bea}{\begin{eqnarray}\label}
\newcommand{\eea}{\end{eqnarray}}
\newcommand{\bas}{\begin{eqnarray*}}
\newcommand{\eas}{\end{eqnarray*}}
\newcommand{\bit}{\begin{itemize}}
\newcommand{\eit}{\end{itemize}}
\newcommand{\qed}{\hfill$\Box$ \vskip.2cm}
\newcommand{\nn}{\nonumber}
\newcommand{\R}{\mathbb{R}}
\newcommand{\N}{\mathbb{N}}
\newcommand{\pO}{\partial\Omega}

\newcommand{\eps}{\varepsilon}

\newcommand{\io}{\int_\Omega}

\newcommand{\abs}{\\[5pt]}

\newcommand{\F}{{\cal F}}
\newcommand{\D}{{\cal D}}
\newcommand{\B}{{\cal B}}
\newcommand{\set}{{\cal S}(m,M,B,\kappa)}
\newcommand{\tm}{T_{max}(u_0,v_0)}
\newcommand{\om}{\overline{m}}
\begin{document}
\title{Finite-time blow-up in the higher-dimensional parabolic-parabolic Keller-Segel system}
\author{
\and
Michael Winkler\footnote{michael.winkler@math.uni-paderborn.de}\\
{\small Institut f\"ur Mathematik, Universit\"at Paderborn,}\\
{\small 33098 Paderborn, Germany} }
\date{}
\maketitle
\begin{abstract}
\noindent 
  We study the Neumann initial-boundary value problem for the fully parabolic Keller-Segel system
  \bas
	\left\{ \begin{array}{l}
	u_t=\Delta u - \nabla \cdot (u\nabla v), \qquad x\in\Omega, \ t>0, \\[1mm]
	v_t=\Delta v-v+u, \qquad x\in\Omega, \ t>0,
	\end{array} \right.
  \eas
  where $\Omega$ is a ball in $\R^n$ with $n\ge 3$.\\
  It is proved that for any prescribed $m>0$ there exist radially symmetric positive initial data 
  $(u_0,v_0) \in C^0(\bar\Omega) \times W^{1,\infty}(\Omega)$ with $\io u_0=m$
  such that the corresponding solution blows up in finite time. 
  Moreover, by providing an essentially explicit blow-up criterion 
  it is shown that within the space of all radial functions, the set of such blow-up enforcing initial data 
  indeed is large in an appropriate sense; in particular, this set is dense with respect to
  the topology of $L^p(\Omega) \times W^{1,2}(\Omega)$ for any $p \in (1,\frac{2n}{n+2})$.\abs
\end{abstract}
\newpage
\section{Introduction}\label{intro}
We consider the parabolic initial-boundary value problem
\be{0}
    \left\{ \begin{array}{l}
	u_t=\Delta u - \nabla \cdot (u\nabla v), \qquad x\in \Omega, \ t>0, \\[1mm]
	v_t=\Delta v-v+u, \qquad x\in\Omega, \ t>0, \\[1mm]
	\frac{\partial u}{\partial\nu}=\frac{\partial v}{\partial\nu}=0, \qquad x\in\pO, \ t>0, \\[1mm]
	u(x,0)=u_0(x), \quad v(x,0)=v_0(x), \qquad x\in\Omega,
    \end{array} \right.
\ee
in a bounded domain $\Omega\subset \R^n$, with given nonnegative initial data $u_0$ and $v_0$.\\
In a celebrated work by Keller and Segel (\cite{keller_segel}), this system was proposed as a macroscopic model for
chemotactic cell migration, that is, for the motion of cells which, besides diffusing randomly, partly orient
their movement towards increasing concentrations of a chemical signal substance.
In this context, $u=u(x,t)$ then denotes the cell density, whereas $v=v(x,t)$ represents the concentration
of the chemical.\abs
In prototypical processes such as the collective behavior of the slime mold {\em Dictyostelium Discoideum},
the signal is produced by the cells themselves, and the possibly most striking consequence thereof appears to be 
the ability of cell populations to spontaneously form aggregates in small spatial regions after a finite time.
Correspondingly, verifying the validity of any mathematical model for such processes is closely linked to
investigating its capability to adequately describe such phenomena of self-organization.
A commonly accepted mathematical concept for this consists of identifying the emergence of aggregation
with the collapse of the corresponding solution into a singularity with respect to the norm
in $L^\infty(\Omega)$ (\cite{horstmann_dmv}, \cite{nanjundiah}). \abs
{\bf The challenge of proving blow-up.}\quad
%
%
%
%
Accordingly, since the work of Keller and Segel 		
extensive mathematical efforts have been undertaken to detect unbounded solutions in (\ref{0})
or, more generally, to determine conditions on the initial data and on certain parameters which either guarantee or rule out
the existence of such blow-up solutions in (\ref{0}) and related models.\\
It turned out that the mathematical difficulties linked to the subtle task of finding unbounded solutions
can significantly be reduced upon replacing (\ref{0}) with associated parabolic-elliptic variants,
the second equation of which being either
\be{e1}
	0=\Delta v - v + u, \qquad x\in\Omega, \ t>0,
\ee
or
\be{e2}
	0=\Delta v - \om + u, \qquad x\in\Omega, \ t>0.
\ee
Here $\om:=\frac{m}{|\Omega|}$, where $m:=\io u_0$ denotes the total mass of cells which remains constant in time in the
sense that $\io u(x,t)dx \equiv m$ for $t>0$.
According to the applicability of much a larger repertoire of mathematical tools, the knowledge
for the resulting simplified systems has been rather complete in respect of the occurrence of blow-up for quite a while.
In fact, it was shown in \cite{jaeger_luckhaus}, \cite{nagai1995} and \cite{nagai2001} that the corresponding
initial-boundary value problems in the spatially two-dimensional setting indeed possess some solutions which blow up
in finite time provided that the mass $m$ is large enough and concentrated around some point to a suitable extent,
whereas if $m$ is small then solutions remain bounded;
the precise threshold values for the mass could be identified
as $8\pi$ in the radially symmetric setting and $4\pi$ in the general case (cf.~\cite{nagai1995}, 
\cite{nagai2001} and \cite{biler1999}, for instance, and also \cite{horstmann_dmv} for a survey). 
In the three-dimensional framework, finite-time blow-up in the systems related to 
(\ref{e1}) and (\ref{e2}) may occur for arbitrarily small values of $m$, meaning that no mass threshold
for aggregation exists in that case (\cite{nagai2000}, \cite{herrero_medina_velazquez}).\\
Actually, the understanding of these simplified systems is even elaborate: For instance, 
for both the two- and three-dimensional cases the studies in
\cite{herrero_velazquez1996_matann} and \cite{herrero_medina_velazquez}
provide examples of unbounded solutions the asymptotic behavior of which can be described rather precisely near their
blow-up time.\abs
In contrast to this, the knowledge on the full parabolic-parabolic system (\ref{0}) appears to be much less
comprehensive except for the case $n=1$ where blow-up is entirely ruled out (\cite{osaki_yagi}).
For example, a counterpart of the above
two-dimensional mass threshold phenomenon could rigorously be proved to exist only in a weakened sense.
Namely, it is known that all solutions of (\ref{0}) with mass satisfying $m<4\pi$ remain bounded for all times
(\cite{gajewski_zacharias}, \cite{nagai_senba_yoshida}), while for any $\eps>0$
there exist some unbounded solutions mith mass $m<4\pi+\eps$ (\cite{horstmann_wang}). 
In general it is not known, however, whether the blow-up time of the latter solutions is finite or infinite; it is thus 
conceivable that these solutions are global in time and become unbounded only in the large time limit. 
Only in the radial symmetric setting certain particular 
solutions of (\ref{0}) satisfying $m > 8 \pi$ were constructed which blow up
in finite time and, moreover, their asymptotic behavior near the blow-up time was described in \cite{herrero_velazquez1997_pisa}.
However, this does not clarify whether or not this phenomenon is exceptional; indeed, nothing is known about the 
size or the structure of the set of initial data enforcing finite-time blow-up.\\
Corresponding confinements in the analogy with the parabolic-elliptic systems appear also in the higher-dimensional
setting: It has been shown in \cite{win_jde} 
that when $n\ge 3$, (\ref{0}) possesses unbounded solutions with arbitrarily small total mass
of cells, but again it has been left open there whether or not the associated blow-up time is finite.
To the best of our knowledge, not even a single example of a solution of (\ref{0}) which undergoes finite-time
blow-up is known in the case $n\ge 3$. \abs
As a more general observation, let us note that
due to their diffusive and hence essentially dissipative structure, systems of 
type (\ref{0}) are amenable to various powerful techniques of regularity theory.
This becomes manifest in the literature not only on the original model (\ref{0}), but
also on quite a large variety of related models involving, for instance, different mechanisms of diffusion and chemotactic
cross-diffusion in the equation for the cell density.
Namely, there is a fairly rich literature addressing boundedness issues 
in such systems of both parabolic-parabolic and parabolic-elliptic type
(see e.g.~\cite{senba_suzuki_bounded}, \cite{kowalczyk}, \cite{kowalczyk_szymanska}, \cite{kozono_sugiyama_jde},
\cite{corrias_perthame}, \cite{calvez_carrillo}, \cite{horstmann_win}, \cite{taowin_jde} and also the survey 
\cite{hillen_painter2009}).
In some special cases even more subtle analytical results on bounded solutions are available, such as e.g.~on
attractors (\cite{wrzosek}) or on two-dimensional
forward self-similar solutions with supercritical mass
(cf.~\cite{biler_corrias_dolbeaut} and the references therein). 
As opposed to this, the few results addressing blow-up also in the setting of nonlinear diffusion 
concentrate on parabolic-elliptic simplifications
(see \cite{cieslak_laurencot2009}, \cite{cieslak_laurencot_dcds}, 
\cite{blanchet_carrillo_laurencot}, \cite{cieslak_win}, \cite{djie_win}, for instance)
or on unboundedness in possibly infinite time (\cite{horstmann_win}, \cite{win_mmas}),
the apparently only exception being a recent
result on finite-time blow-up in a quasilinear one-dimensional Keller-Segel system with sufficiently weak
nonlinear diffusion of cells and sufficiently fast diffusion of chemoattractant (\cite{cieslak_laurencot_ann_ihp}).\abs
{\bf Main results.} \quad
In view of the underlying biological background, we find it worthwhile to firstly investigate whether
cell aggregation in the mathematically extreme flavor of finite-time blow-up 
at all occurs in (\ref{0}) also in the three-dimensional setting, and secondly to make sure whether this
is a rare phenomenon which can only be expected 
under very special assumptions on the initial framework, or whether the set of initial data enforcing
a finite-time collapse is rich in an appropriate sense.
Accordingly, the purpose of the present paper consists of deriving sufficient conditions on $(u_0,v_0)$ which
lead to finite-time blow-up of solutions to (\ref{0}) in the more general case $n\ge 3$.\\
In order to formulate our main results in this direction, and thereby demonstrate that moreover
these criteria will be essentially explicit, let us recall
that any solution of (\ref{0}) satisfies the energy inequality
\be{energy}
	\frac{d}{dt} \F(u(\cdot,t),v(\cdot,t)) \le - \D (u(\cdot,t),v(\cdot,t)) 
	\qquad \mbox{for all } t \in (0,\tm)
\ee
where $\tm \in (0,\infty]$ denotes the maximum existence time of $(u,v)$ and where,
for arbitrary smooth positive functions $u$ and $v$, the energy is defined by
\be{F}
	\F(u,v):=\frac{1}{2} \io |\nabla v|^2 + \frac{1}{2} \io v^2 - \io uv + \io u\ln u
\ee
and the dissipation rate is given by
\be{D}
	\D(u,v):=\io v_t^2 + \io u \cdot \Big| \frac{\nabla u}{u} - \nabla v \Big|^2
\ee
(see \cite{nagai_senba_yoshida}, \cite{win_jde} and also Lemma \ref{lem_loc_exist} below).\abs
This energy inequality plays an essential role in deriving boundedness of solutions in subcritical cases
(see \cite{nagai_senba_yoshida}, for instance). But also the detection of unbounded (possibly global) solutions
in \cite{horstmann_wang} and \cite{win_jde} crucially relies on the use of (\ref{energy}) 
through an indirect argument: Indeed, if $\F(u_0,v_0)$ is a sufficiently large negative number then $(u,v)$ cannot be 
both global and bounded, because due to (\ref{energy}) any such solution must approach a low-energy equilibrium
which ruled out by a corresponding a-priori estimate for energies of steady states. 
This type of reasoning, well-established in scalar parabolic problems (\cite{levine_survey}) and also applicable 
in larger classes
of quasilinear Keller-Segel models (\cite{horstmann_wang}, \cite{horstmann_win}, \cite{win_mmas}), 
evidently cannot give any information on the actual blow-up time,
and hence cannot rule out the possibility that $\tm=\infty$.
However, it will turn out in this work that a more subtle analysis of (\ref{energy}) 
can be used to derive a basically explicit sufficient condition on the initial data which ensure that
finite-time blow-up occurs.
More precisely, the first of our main results reads as follows.
\begin{theo}\label{theo12}
  Let $\Omega=B_R \subset \R^n$ with some $n\ge 3$ and $R>0$, and let $m>0$ and $A>0$. 
  Then there exist $T(m,A)>0$ and $K(m,A)>0$ with the property that given any
  $(u_0,v_0)$ from the set
  \bea{12.1}
	\B (m,A) &:=& \bigg\{
	(u_0,v_0) \in C^0(\bar\Omega) \times W^{1,\infty}(\Omega) \ \bigg| \
	\mbox{$u_0$ and $v_0$ are radially symmetric and positive in $\bar\Omega$} \nn\\
	& & \hspace*{20mm}
	\mbox{with $\io u_0=m$, $\|v_0\|_{W^{1,2}(\Omega)} \le A$ and $\F(u_0,v_0) \le -K(m,A)$} \bigg\},
  \eea
  for the corresponding solution $(u,v)$ of (\ref{0}) we have $\tm \le T(m,A) <\infty$; that is, $(u,v)$ blows up before 
  or at time $T(m,A)$.
\end{theo}
Secondly, we shall address the question in how far the above set of low-energy initial data can be considered large.
In fact, this set turns out to be even dense in the space of positive
functions in $C^0(\bar\Omega) \times W^{1,\infty}(\Omega)$ when equipped with an appropriate topology:
\begin{theo}\label{theo15}
  Let $\Omega$ be as in Theorem \ref{theo12}, and suppose that $p \in (1,\frac{2n}{n+2})$. 
  Then for each $m>0$ and $A>0$, the set $\B(m,A)$ defined in (\ref{12.1})
  is dense in the space of all radially symmetric positive functions in $C^0(\bar\Omega) \times W^{1,\infty}(\Omega)$
  with respect to the topology in $L^p(\Omega) \times W^{1,2}(\Omega)$. \\
  In particular, for any positive radial $(u_0,v_0) \in C^0(\bar\Omega) \times W^{1,\infty}(\Omega)$ and any 
  $\eps>0$ one can find some radial positive $(u_{0\eps},v_{0\eps}) \in C^0(\bar\Omega) \times W^{1,\infty}(\Omega)$
  such that 
  \bas
	\|u_{0\eps} - u_0\|_{L^p(\Omega)} + \|v_{0\eps}-v_0\|_{W^{1,2}(\Omega)} < \eps,
  \eas
  but such that the solution $(u_\eps,v_\eps)$ of (\ref{0}) with initial data $(u_\eps,v_\eps)|_{t=0}=(u_{0\eps},v_{0\eps})$
  blows up in finite time.
\end{theo}
Let us underline that to the best of our knowledge, this is the first result asserting the occurrence of
finite-time blow-up in the Keller-Segel system (\ref{0}) in space dimension $n\ge 3$.
But Theorem \ref{theo12} and Theorem \ref{theo15} evidently go much further:
For instance, they especially say that each of the constant steady states $(u,v) \equiv (c,c)$, $c>0$, is highly 
unstable in that any of its neighborhoods in the above topology contains initial data which evolve into a singularity
in finite time.\abs
{\bf Plan of the paper.} \quad
%
%
Our technical approach is based on the idea to estimate the dissipated quantity in (\ref{energy}) from below
in order to turn (\ref{energy}) into an inequality of the form
\bas
	\frac{d}{dt} \Big(-\F(u(\cdot,t),v(\cdot,t)) \Big)
	\ge \Big( c \cdot (-\F(u(\cdot,t),v(\cdot,t))) -1 \Big)_+^\lambda \qquad \mbox{for all } t\in (0,\tm)
\eas
with some $\lambda>1$ and $c>0$.
We shall thus be concerned with deriving an upper bound for $\F(u,v)$ in terms of a sublinear power of $\D(u,v)$,
and the main step towards this will be provided by the estimate
\bas
	\io uv \le C \cdot 
	\bigg( \Big\|\Delta v-v+u\Big\|_{L^2(\Omega)}^{2\theta} 
	+\Big\|\frac{\nabla u}{\sqrt{u}}-\sqrt{u}\nabla v\Big\|_{L^2(\Omega)} +1 \bigg)
\eas
with some $\theta \in (0,1)$ and $C>0$, to be given in Lemma \ref{lem8}. 
Here, $u$ and $v$ will be allowed to be rather arbitrary smooth positive
radial functions satisfying appropriate mass constraints and an additional pointwise upper estimate for $v$
that can be shown to be fulfilled by 
the component $v(\cdot,t)$ for any $t\in (0,\tm)$ of any solution of (\ref{0}) in question (Corollary \ref{cor3}).\\
The application to the parabolic problem is then straightforward (see Section \ref{sect_blowup}),
whereas Theorem \ref{theo15} will be proved using an explicit construction of appropriate 
initial data in Section \ref{sect_dense}.
\mysection{Preliminaries}
To begin with, let us collect some basic statements on local well-posedness and elementary properties of solutions
to (\ref{0}).
\begin{lem}\label{lem_loc_exist}
  Let $(u_0,v_0) \in C^0(\bar\Omega) \times W^{1,\infty}(\Omega)$ be radially symmetric and positive in $\bar\Omega$,
  and fix $q\in (n,\infty)$.
  Then there exist $\tm \in (0,\infty]$ and a uniquely determined couple $(u,v)$ of radially symmetric functions,
  satisfying the inclusions
  \bas
	& & u \in C^0([0,\tm);C^0(\bar\Omega)) \cap C^{2,1}(\bar\Omega \times (0,\tm)) \qquad \mbox{and} \\
	& & v \in C^0([0,\tm);W^{1,q}(\Omega)) \cap C^{2,1}(\bar\Omega \times (0,\tm)),
  \eas
  which solves (\ref{0}) classically in $\Omega \times (0,\tm)$ and has the extensibility 
  property
  \bas
	\mbox{either $\tm=\infty$, or } \|u(\cdot,t)\|_{L^\infty(\Omega)} 		
	\to \infty \qquad \mbox{as } t\nearrow \tm.
  \eas
  In addition, this solution fulfils
  \be{mass}
	\io u(x,t)dx = \io u_0 \qquad \mbox{for all } t\in (0,\tm)
  \ee
  and
  \be{mass_v}
	\io v(x,t)dx \le \max \Big\{ \io u_0, \io v_0 \Big\}
	\qquad \mbox{for all } t\in (0,\tm),
  \ee
  and moreover the energy inequality (\ref{energy}) holds.
\end{lem}
\proof 
  The statements concerning existence, uniqueness, regularity and extensibility are well-known, 
  and thus for details covering the present and more general frameworks 
  we may refer the reader to \cite{horstmann_win}, \cite{biler1999} and \cite{yagi1997}, for instance.\\
  The identity (\ref{mass}) immediately follows from integration of the first equation in (\ref{0}) in space,
  whereupon integrating the second one yields
  \bas
	\frac{d}{dt} \io v(x,t)dx = - \io v(x,t)dx + \io u_0 \qquad \mbox{for all } t\in (0,\tm).
  \eas
  Combined with a straightformward ODE comparison, this proves (\ref{mass_v}).
\qed
The following consequences of the Gagliardo-Nirenberg inequality and Young's inequality are immediate.
Since they will be used in several places in the sequel, let us briefly state them separately and refer to
\cite{friedman_book} for the underlying interpolation estimates.
\begin{lem}\label{lem13}
  There exists $C>0$ such that
  \be{13.1}
	\|\varphi\|_{L^2(\Omega)} \le C\|\nabla \varphi\|_{L^2(\Omega)}^\frac{n}{n+2} \|\varphi\|_{L^1(\Omega)}^\frac{2}{n+2}
	+ C\|\varphi\|_{L^1(\Omega)}
	\qquad \mbox{for all } \varphi \in W^{1,2}(\Omega).
  \ee
  Moreover, for each $\eps>0$ one can find $C(\eps)>0$ with the property that
  \be{13.2}
	\|\varphi\|_{L^2(\Omega)}^2 \le \eps \|\nabla \varphi\|_{L^2(\Omega)}^2 + C(\eps) \|\varphi\|_{L^1(\Omega)}^2
	\qquad \mbox{for all } \varphi \in W^{1,2}(\Omega).
  \ee
\end{lem}
\mysection{A pointwise upper bound for solutions of (\ref{0})}
Let us first adapt a basically well-known regularity property of the second solution component $v$ which
is a straightforward consequence of standard parabolic regularity arguments and thereby it
does in fact not require any symmetry assumption on the initial data.
\begin{lem}\label{lem1}
  Let $p \in (1,\frac{n}{n-1})$. Then there exists $C(p)>0$ such that for any choice of positive functions
  $u_0 \in C^0(\bar\Omega)$ and $v_0 \in W^{1,\infty}(\Omega)$, the solution of (\ref{0}) satisfies
  \be{1.1}
	\|\nabla v(\cdot,t)\|_{L^p(\Omega)} \le C(p) \cdot \Big( \|u_0\|_{L^1(\Omega)} 
	+ \|\nabla v_0\|_{L^2(\Omega)} \Big)
	\qquad \mbox{for all } t \in (0,\tm).
  \ee
\end{lem}
\proof
  It s well-known (cf.~\cite{win_jde}, for instance) that the Neumann heat semigroup $(e^{t\Delta})_{t\ge 0}$
  in $\Omega$ has the property
  \be{1.2}
	\|\nabla e^{t\Delta}\varphi\|_{L^p(\Omega)} \le c_1 t^{-\frac{1}{2}-\frac{n}{2}(1-\frac{1}{p})} 
	\|\varphi\|_{L^1(\Omega)}
	\qquad \mbox{for all } \varphi \in L^1(\Omega)
  \ee
  with some $c_1>0$.
  Moreover, since it can easily be checked that $\frac{d}{dt} \|\nabla e^{t\Delta} \varphi\|_{L^2(\Omega)}^2 \le 0$, using
  the H\"older inequality along with the fact that $p<2$ we can find $c_2>0$ such that
  \be{1.3}
	\|\nabla e^{t\Delta}\varphi\|_{L^p(\Omega)} 
	\le c_2 \|\nabla e^{t\Delta}\varphi\|_{L^2(\Omega)} \le c_2 \|\nabla \varphi\|_{L^2(\Omega)}
	\qquad \mbox{for all } \varphi \in W^{1,2}(\Omega).
  \ee
  Applying (\ref{1.2}) and (\ref{1.3}) to the variation-of-constants representation of $v$,
  \bas
	v(\cdot,t)=e^{t(\Delta-1)} v_0 + \int_0^t e^{(t-s)(\Delta-1)} u(\cdot,s) ds,
	\qquad t\in (0,\tm),
  \eas
  we obtain
  \bas
	\|\nabla v(\cdot,t)\|_{L^p(\Omega)} 
	&\le& c_2 \|\nabla v_0\|_{L^2(\Omega)}
	+ c_1 \int_0^t (t-s)^{-\frac{1}{2}-\frac{n}{2}(1-\frac{1}{p})} \cdot e^{-(t-s)} \|u(\cdot,s)\|_{L^1(\Omega)} ds\\
	&=& c_2 \|\nabla v_0\|_{L^2(\Omega)}
	+ c_1 \|u_0\|_{L^1(\Omega)} \cdot \int_0^t \sigma^{-\frac{1}{2}-\frac{n}{2}(1-\frac{1}{p})} e^{-\sigma} d\sigma \\
	 \qquad \mbox{for all } t\in (0,\tm),
  \eas	
  because $\|u(\cdot,t)\|_{L^1(\Omega)} = \|u_0\|_{L^1(\Omega)}$ for all $t\in (0,\tm)$ by (\ref{mass}).
  Now since our restriction $p<\frac{n}{n-1}$ ensures that $\frac{1}{2}+\frac{n}{2}(1-\frac{1}{p})<1$, this implies
  \bas
	\|\nabla v(\cdot,t)\|_{L^p(\Omega)} 
	&\le& c_2 \|\nabla v_0\|_{L^2(\Omega)}
	+ c_1 \|u_0\|_{L^1(\Omega)} \cdot \Gamma \Big(\frac{1}{2}-\frac{n}{2}\Big(1-\frac{1}{p}\Big) \Big)
	\qquad \mbox{for all } t\in (0,\tm)
  \eas	
  and thereby proves (\ref{1.1}).
\qed
As a consequence of Lemma \ref{lem1} and (\ref{mass_v}), 
in the case when $(u,v)$ is radially symmetric we obtain a pointwise upper bound for $v$
which is valid up to the blow-up time and hence gives a first, though rather rough, information on what
might finally be called the spatial blow-up profile of $(u,v)$.
\begin{lem}\label{lem2}
  Let $p\in (1,\frac{n}{n-1})$. Then there exists $C(p)>0$ such that
  whenever $u_0 \in C^0(\bar\Omega)$ and $v_0 \in W^{1,\infty}(\Omega)$ are positive and radially symmetric, 
  the solution of (\ref{0}) satisfies
  \be{2.1}
	v(r,t) \le C(p) \cdot \Big( \|u_0\|_{L^1(\Omega)} + \|v_0\|_{L^1(\Omega)}
	+ \|\nabla v_0\|_{L^2(\Omega)} \Big) \cdot r^{-\frac{n-p}{p}}
	\qquad \mbox{for all } (r,t) \in (0,R) \times (0,\tm).
  \ee
\end{lem}
\proof
  Abbreviating $M:=\max \Big\{\io u_0, \io v_0 \Big\}$, from (\ref{mass_v}) we know that
  $\io v(\cdot,t) \le M$ for all $t\in (0,\tm)$. Thus, for each $t\in (0,\tm)$ we can pick some
  $r_0(t) \in (\frac{R}{2},R)$ such that 
  \bas
	v(r_0(t),t) \le \frac{M}{|B_R \setminus B_\frac{R}{2}|},
  \eas
  for supposing the opposite would lead to the absurd conclusion
  \bas
	\io v(\cdot,t) \ge \int_{B_R \setminus B_\frac{R}{2}} v(\cdot,t) 
	> \int_{B_R \setminus B_\frac{R}{2}} \frac{M}{|B_R \setminus B_\frac{R}{2}|} =M.
  \eas
  Therefore, using the H\"older inequality and Lemma \ref{lem1}, we can find $c_1(p)>0$ such that
  \bea{2.2}
	v(r,t)-v(r_0(t),t)
	&=& \int_{r_0(t)}^r v_r(\rho,t) d\rho \nn\\
	&\le& \bigg| \int_{r_0(t)}^r \rho^{n-1} |v_r(\rho,t)|^p d\rho \bigg|^\frac{1}{p}
	\cdot \bigg| \int_{r_0(t)}^r \rho^{-\frac{n-1}{p-1}} d\rho \bigg|^\frac{p-1}{p} \nn\\
	&\le& c_1(p) \cdot (\|u_0\|_{L^1(\Omega)} + \|\nabla v_0\|_{L^2(\Omega)}) 
	\cdot \bigg| \int_{r_0(t)}^r \rho^{-\frac{n-1}{p-1}} d\rho \bigg|^\frac{p-1}{p}
  \eea
  for all $(r,t) \in (0,R) \times (0,\tm)$. 
  Now since $p<\frac{n}{n-1}<n$, for small $r$ we can estimate
  \bea{2.3}
	\bigg| \int_{r_0(t)}^r \rho^{-\frac{n-1}{p-1}} d\rho \bigg|^\frac{p-1}{p}
	&=& \bigg( \frac{r^{-\frac{n-p}{p-1}} - (r_0(t))^{-\frac{n-p}{p-1}}}{\frac{n-p}{p-1}} \bigg)^\frac{p-1}{p} \nn\\
	&\le& \Big(\frac{p-1}{n-p}\Big)^\frac{p-1}{p} \cdot r^{-\frac{n-p}{p}}
	\qquad \mbox{for all } r\in (0,r_0(t)],
  \eea
  whereas for large $r$	we similarly find that
  \bea{2.4}
	\bigg| \int_{r_0(t)}^r \rho^{-\frac{n-1}{p-1}} d\rho \bigg|^\frac{p-1}{p}
	&\le& \Big(\frac{p-1}{n-p}\Big)^\frac{p-1}{p} \cdot (r_0(t))^{-\frac{n-p}{p}}
	\le \Big(\frac{p-1}{n-p}\Big)^\frac{p-1}{p} \cdot \Big(\frac{R}{2}\Big)^{-\frac{n-p}{p}} \nn\\
	&\le& \Big(\frac{p-1}{n-p}\Big)^\frac{p-1}{p} \cdot 2^\frac{n-p}{p} r^{-\frac{n-p}{p}}
	\qquad \mbox{for all } r\in (r_0(t),R),
  \eea
  because $r_0(t)>\frac{R}{2}$. As finally
  \bas
	v(r_0(t),t) \le \frac{M}{|B_R \setminus B_\frac{R}{2}|} 
	\le \frac{M}{|B_R \setminus B_\frac{R}{2}|} \cdot R^\frac{n-p}{p} r^{-\frac{n-p}{p}}
	\qquad \mbox{for all } r\in (0,R),
  \eas
  (\ref{2.2})-(\ref{2.4}) imply (\ref{2.1}).
\qed 
Adjusting $p$ in Lemma \ref{lem2} appropriately, we can achieve an estimate showing that the singularity
of $v$ can essentially not be stronger than that of the fundamental solution of the Laplacian.
\begin{cor}\label{cor3}
  Let $\kappa>n-2$. Then one can find $C(\kappa)>0$ such that for all radially symmetric and positive functions
  $u_0 \in C^0(\bar\Omega)$ and $v_0\in W^{1,\infty}(\Omega)$, the corresponding solution of (\ref{0}) satisfies
  \be{3.1}
	v(r,t) \le C(\kappa) \cdot \Big(\|u_0\|_{L^1(\Omega)} + \|v_0\|_{L^1(\Omega)} + \|\nabla v_0\|_{L^2(\Omega)} \Big)
	\cdot r^{-\kappa}
	\qquad \mbox{for all } (r,t) \in (0,R) \times (0,\tm).
  \ee
\end{cor}
\proof 
  Since $\kappa>n-2$, we have $\frac{n}{\kappa+1}<\frac{n}{n-1}$, so that it is possible to fix $p>1$ such that
  $\frac{n}{\kappa+1} \le p < \frac{n}{n-1}$. An application of Lemma \ref{lem2} then easily yields (\ref{3.1}),
  because $p\ge \frac{n}{\kappa+1}$ implies $\frac{n-p}{p} \le \kappa$.
\qed
\mysection{An estimate for $\io uv$ in terms of the dissipation rate}\label{sect_estimates}
Guided by our knowledge on the solutions of (\ref{0}) gained above, 
in asserting a lower estimate of the desired form
\bas
	\frac{\F(u,v)}{\D^\theta(u,v)+1} \ge - C
\eas
with some $\theta \in (0,1)$ and $C>0$,
we shall concentrate henceforth on positive radial functions satisfying
the mass constraints
\be{m}
	\io u = m \qquad \mbox{and} \qquad
	\io v \le M
\ee
and the additional pointwise restriction
\be{B}
	v(x) \le B|x|^{-\kappa} \qquad \mbox{for all } x\in\Omega,
\ee
where $m>0, M>0, B>0$ and $\kappa>n-2$ are given fixed parameters.\abs
More precisely, our goal will be to derive an inequality of the form
\be{bound}
	\frac{\F(u,v)}{\D^\theta(u,v)+1} \ge - C(m,M,B,\kappa)
	\qquad \mbox{for all } (u,v) \in \set
\ee
with some $\theta \in (0,1)$ and $C(m,M,B,\kappa)>0$, where
\bea{S}
	\set &:=& \bigg\{ (u,v) \in C^1(\bar\Omega) \times C^2(\bar\Omega) \ \bigg| \
	\mbox{$u$ and $v$ are positive and radially symmetric} \nn\\
	& & \hspace*{24mm}
	\mbox{with $\frac{\partial v}{\partial\nu}=0$ on $\pO$ and such that (\ref{m}) and (\ref{B}) hold} \bigg\}
\eea
and $\F$ and $\D$ are as defined in (\ref{F}) and (\ref{D}), respectively.\\
In view of the latter,		
establishing (\ref{bound}) essentially amounts to 
showing that the integral $\io uv$ is bounded from above by a sublinear power of 
the sum of the norms in $L^2(\Omega)$ of 
the functions $f$ and $g$ which for convenience in subsequent notation are introduced by abbreviating
\be{f}
	f:=-\Delta v + v - u
\ee
and
\be{g}
	g:=\Big(\frac{\nabla u}{\sqrt{u}}-\sqrt{u}\nabla v \Big) \cdot \frac{x}{|x|} \qquad (\mbox{for } x\ne 0)
\ee
for $(u,v)\in \set$; since $(u,v)$ is radial, these definitions actually reduces to the identities 
$f=-r^{1-n}(v^{n-1}v_r)_r + v - u$ and $g=\frac{u_r}{\sqrt{u}} - \sqrt{u}v_r$.\\
The essential step toward (\ref{bound}) will be contained in the following main result of this section.
\begin{lem}\label{lem8}
  There exists $C(m,M,B,\kappa)>0$ such that for all $(u,v)\in\set$ we have
  \be{8.1}
	\io uv \le C(m,M,B,\kappa) \cdot 
	\bigg( \Big\|\Delta v-v+u\Big\|_{L^2(\Omega)}^{2\theta} 
	+ \Big\|\frac{\nabla u}{\sqrt{u}}-\sqrt{u}\nabla v\Big\|_{L^2(\Omega)} +1 \bigg)
  \ee
  with $\theta \in (\frac{1}{2},1)$ given by (\ref{6.2}).
\end{lem}
%
%
%
%
%
%
%
%
The proof of Lemma \ref{lem8} will be accomplished through a series of auxiliary statements.
The first of these shows that proving (\ref{8.1}) actually amounts to estimating $\io |\nabla v|^2$.
\begin{lem}\label{lem7}
  There exists $C(M)>0$ such that for all $(u,v) \in \set$ we have
  \be{7.1}
	\io uv \le 2 \io |\nabla v|^2 + C(M) \cdot \bigg( 
	\Big\| \Delta v - v + u \Big\|_{L^2(\Omega)}^\frac{2n+4}{n+4} +1 \bigg).
  \ee
\end{lem}
\proof
  We use the notation in (\ref{f}) and multiply the latter by $v$ to obtain upon integrating by parts over $\Omega$ that
  \be{7.2}
	\io uv = \io |\nabla v|^2 + \io v^2 - \io fv.
  \ee
  In order to estimate the right-hand side appropriately, we note that by Lemma \ref{lem13} and (\ref{m}) there exists
  $c_1=c_1(M)>0$ and $c_2=c_2(M)>0$ such that
  \be{7.3}
	\|v\|_{L^2(\Omega)} \le c_1 \Big(\|\nabla v\|_{L^2(\Omega)}^\frac{n}{n+2}+1 \Big)
  \ee
  and
  \be{7.4}
	\io v^2 \le \frac{1}{2} \io |\nabla v|^2 + c_2.
  \ee
  Furthermore, combining (\ref{7.3}) with the Cauchy-Schwarz inequality and Young's inequality applied
  with exponents $\frac{2n+4}{n}$ and $\frac{2n+4}{n+4}$ provides $c_3=c_3(M)>0$ such that
  \bea{7.5}
	- \io fv
	&\le& \|f\|_{L^2(\Omega)} \|v\|_{L^2(\Omega)} \nn\\
	&\le& c_1 \cdot (\|\nabla v\|_{L^2(\Omega)}^\frac{n}{n+2}+1) \cdot \|f\|_{L^2(\Omega)} \nn\\
	&\le& \frac{1}{2} \io |\nabla v|^2 
	+ c_3 \|f\|_{L^2(\Omega)}^\frac{2n+4}{n+4} + c_1 \|f\|_{L^2(\Omega)}.
  \eea
  Since $\frac{2n+4}{n+4}>1$, again by Young's inequality we can find $c_4=c_4(M)$ fulfilling
  \bas
	c_1 \|f\|_{L^2(\Omega)} \le \|f\|_{L^2(\Omega)}^\frac{2n+4}{n+4} + c_4,
  \eas
  whereupon it becomes clear that (\ref{7.2}), (\ref{7.4}) and (\ref{7.5}) imply (\ref{7.1}).
\qed
Accordingly, our next goal is to bound $\io |\nabla v|^2$ appropriately. 
This will be done by splitting this expression into an integral over a small inner ball $B_{r_0}$ and
a corresponding outer annulus, the precise value of $r_0$ remaining at our disposal until it will be fixed
in Lemma \ref{lem6} below. Let us first concentrate on the outer region.
\begin{lem}\label{lem4}
  Let $r_0\in (0,R)$ and $\eps \in (0,1)$. Then one can find a constant 
  $C(\eps,m,M,B,\kappa)>0$ such that for all $(u,v)\in\set$ the estimate
  \bea{4.1}
	\int_{\Omega \setminus B_{r_0}} |\nabla v|^2 
	\le \eps \io uv + \eps \io |\nabla v|^2 
	+ C(\eps,m,M,B,\kappa) \cdot \bigg\{
	r_0^{-\frac{2n+4}{n}\kappa}
	+ \Big\|\Delta v-v+u\Big\|_{L^2(\Omega)}^\frac{2n+4}{n+4} \bigg\}
  \eea
  holds.
\end{lem}
\proof
  We fix an arbitrary $\alpha \in (0,1)$.
  Then observing that $v\ge 0$, we may multiply the identity (\ref{f}) defining $f$ by $v^\alpha$ and integrate by parts
  over $\Omega$ to obtain
  \be{4.2}
	\alpha \io v^{\alpha-1}|\nabla v|^2 + \io v^{\alpha+1}
	= \io uv^\alpha + \io fv^\alpha.
  \ee
  Here we apply (\ref{B}) and use the fact that $\alpha \in (0,1)$ to estimate
  \bas
	\alpha \io v^{\alpha-1} |\nabla v|^2 \ge \alpha B^{\alpha-1} r_0^{(1-\alpha)\kappa} \cdot
	\int_{\Omega \setminus B_{r_0}} |\nabla v|^2 
  \eas
  and thus infer from (\ref{4.2}) upon dropping a nonnegative term that
  \bea{4.3}
	\int_{\Omega \setminus B_{r_0}} |\nabla v|^2 
	&\le& \frac{B^{1-\alpha}}{\alpha} r_0^{-(1-\alpha)\kappa} \io uv^\alpha \,
	+ \, \frac{B^{1-\alpha}}{\alpha} r_0^{-(1-\alpha)\kappa} \io fv^\alpha.
  \eea
  Now according to Young's inequality, to each $\eta>0$ there corresponds some $c_1(\eta,B)>0$ such that
  \be{4.4}
	\frac{B^{1-\alpha}}{\alpha} r_0^{-(1-\alpha)\kappa} v^\alpha(r)
	\le \eta v(r) + c_1(\eta,B) r_0^{-\kappa}
	\qquad \mbox{for all } r\in (0,R),
  \ee
  which applied to $\eta:=\eps$ yields
  \bea{4.44}
	\frac{B^{1-\alpha}}{\alpha} r_0^{-(1-\alpha)\kappa} \io uv^\alpha
	&\le& \eps \io uv 
	+ c_1(\eps,B) r_0^{-\kappa} \io u \nn\\
	 &=& \eps \io uv 
	+ c_1(\eps,B) m r_0^{-\kappa} \nn\\
	&\le& \eps \io uv 
	+ c_1(\eps,B) m R^{\frac{n+4}{n}\kappa} r_0^{-\frac{2n+4}{n}\kappa}
  \eea
  in view of the nonnegativity of $u$ and (\ref{m}).\\
  Moreover, an application of (\ref{4.4}) to $\eta:=1$ shows that
  \be{4.5}
	\frac{B^{1-\alpha}}{\alpha} r_0^{-(1-\alpha)\kappa} \io fv^\alpha
	\le \io |f|v + c_1(1,B) r_0^{-\kappa} \io |f|,
  \ee
  where by the Cauchy-Schwarz inequality we have
  \bas
	\io |f|v \le \|f\|_{L^2(\Omega)} \|v\|_{L^2(\Omega)}
	\qquad \mbox{and} \qquad
	\io |f| \le \sqrt{|\Omega|} \|f\|_{L^2(\Omega)}.
  \eas
  In order to further estimate the first expression, we invoke Lemma \ref{lem13} which in conjunction with
  (\ref{m}) provides $c_2(M)>0$ such that
  \bas
	\|v\|_{L^2(\Omega)} \le c_2(M) \cdot \Big( \|\nabla v\|_{L^2(\Omega)}^\frac{n}{n+2} + 1 \Big)
    	\le c_2(M) \cdot \Big( \|\nabla v\|_{L^2(\Omega)}^\frac{n}{n+2} + R^\kappa r_0^{-\kappa} \Big),
  \eas
  whence (\ref{4.5}) becomes
  \bas
	\frac{B^{1-\alpha}}{\alpha} r_0^{-(1-\alpha)\kappa} \io fv^\alpha
	\le c_3(M,B,\kappa) \cdot \Big(\|f\|_{L^2(\Omega)} \|\nabla v\|_{L^2(\Omega)}^\frac{n}{n+2}
	+r_0^{-\kappa} \|f\|_{L^2(\Omega)} \Big)
  \eas
  with some $c_3(M,B,\kappa)>0$.
  Here by means of Young's inequality, we can find $c_4(\eps,M,B,\kappa)>0$ and $c_5(M,B,\kappa)>0$ such that
  \bas
	c_3(M,B,\kappa) \|f\|_{L^2(\Omega)} \|\nabla v\|_{L^2(\Omega)}^\frac{n}{n+2}
	\le \eps \|\nabla v\|_{L^2(\Omega)}^2
	+ c_4(\eps,M,B,\kappa) \|f\|_{L^2(\Omega)}^\frac{2n+4}{n+4}
  \eas
  and
  \bas
	c_3(M,B,\kappa) r_0^{-\kappa} \|f\|_{L^2(\Omega)}
	\le \|f\|_{L^2(\Omega)}^\frac{2n+4}{n+4}
	+ c_5(M,B,\kappa) r_0^{-\frac{2n+4}{n}\kappa}.
  \eas
  Therefore, (\ref{4.5}) all in all becomes
  \bas
	\frac{B^{1-\alpha}}{\alpha} r_0^{-(1-\alpha)\kappa} \io fv^\alpha
	\le \eps \io |\nabla v|^2
	+ \Big( c_4(\eps,M,B,\kappa) +1\Big) \cdot \|f\|_{L^2(\Omega)}^\frac{2n+4}{n+4}
	+ c_5(M,B,\kappa) r_0^{-\frac{2n+4}{n}\kappa},
  \eas
  which combined with (\ref{4.44}) and (\ref{4.3}) yields (\ref{4.1}).
\qed
We next estimate $\nabla v$ in the corresponding interior part, where we emphasize the importance of the
factor $r_0$ in the term $r_0  \cdot \|\Delta v-v+u\|_{L^2(\Omega)}^2$ on the right-hand side 
of (\ref{5.1}). Indeed, $r_0$ will eventually be chosen 
in dependence of $\|\Delta v-v+u\|_{L^2(\Omega)}$ in such a way that the above product essentially becomes 
a suitable {\em sub}quadratic power of $\|\Delta v-v+u\|_{L^2(\Omega)}$ (see Lemma \ref{lem6}).
Let us also mention that as (\ref{5.12}) will show, our assumption $n\ge 3$ is crucially needed here.
\begin{lem}\label{lem5}
  There exists $C(m)>0$ such that for any $r_0\in (0,R)$ and all $(u,v)\in\set$ we have
  \be{5.1}
	\int_{B_{r_0}} |\nabla v|^2
	\le C(m) \cdot \bigg\{
	r_0 \cdot \Big\|\Delta v-v+u \Big\|_{L^2(\Omega)}^2
	+\Big\|\frac{\nabla u}{\sqrt{u}}-\sqrt{u}\nabla v\Big\|_{L^2(\Omega)} 
	+\|v\|_{L^2(\Omega)}^2
	+1 \bigg\}.
  \ee
\end{lem}
\proof
  Abbreviating as in (\ref{f}), we rewrite $-\Delta v+v=u+f$ in polar coordinates to see that
  \bas
	(r^{n-1}v_r)_r = -r^{n-1}u - r^{n-1}f + r^{n-1}v, \qquad r\in (0,R),
  \eas
  which we multiply by $r^{n-1}v_r$ to obtain
  \be{5.2}
	\frac{1}{2} \Big((r^{n-1}v_r)^2\Big)_r
	= -r^{2n-2}uv_r - r^{2n-2} fv_r + r^{2n-2} vv_r,
	\qquad r\in (0,R).
  \ee
  In the first term on the right, referring to the notation in (\ref{g}) 
%

  we substitute $v_r=\frac{u_r}{u}-\frac{g}{\sqrt{u}}$, so that
  \be{5.3}
	-r^{2n-2} uv_r=-r^{2n-2}u_r + r^{2n-2} \sqrt{u}g, 
	\qquad r\in (0,R),
  \ee
  whereas for the last term in (\ref{5.2}) we clearly have
  \be{5.4}
	r^{2n-2} vv_r = \frac{1}{2} r^{2n-2} (v^2)_r,
	\qquad r\in (0,R).
  \ee
  As for the expression involving $f$, we pick any $\delta \in (0,\frac{2n-2}{R}]$ and apply Young's inequality to obtain
  \be{5.5}
	-r^{2n-2} fv_r
	\le \frac{\delta}{2} (r^{n-1} v_r)^2 + \frac{1}{2\delta} r^{2n-2} f^2,
	\qquad r\in (0,R).
  \ee
  In light of (\ref{5.3})-(\ref{5.5}), (\ref{5.2}) shows that $y(r):=(r^{n-1}v_r(r))^2, \, r\in [0,R]$, satisfies
  \bas
	y_r \le - 2r^{2n-2} u_r + 2r^{2n-2} \sqrt{u}g
	+ \delta y + \frac{1}{\delta} r^{2n-2} f^2
	+ r^{2n-2} (v^2)_r
	\qquad \mbox{for all } r\in (0,R).
  \eas
  Since $y(0)=0$ thanks to the smoothness of $v$, an integration of this ODI yields
  \bea{5.6}
	r^{2n-2} v_r^2(r) = y(r) &\le&
	-2\int_0^r e^{\delta (r-\rho)} \rho^{2n-2} u_r(\rho) d\rho
	+2\int_0^r e^{\delta (r-\rho)} \rho^{2n-2} \sqrt{u(\rho)} g(\rho) d\rho \nn\\
	& & +\frac{1}{\delta} \int_0^r e^{\delta (r-\rho)} \rho^{2n-2} f^2(\rho) d\rho
	+\int_0^r e^{\delta (r-\rho)} \rho^{2n-2} (v^2)_r(\rho) d\rho
  \eea
  for all $r\in (0,R)$.\\
  Here an integration by parts gives
  \bea{5.7}
	-2\int_0^r e^{\delta (r-\rho)} \rho^{2n-2} u_r(\rho) d\rho
	&=& 4(n-1) \int_0^r e^{\delta(r-\rho)} \rho^{2n-3} u(\rho) d\rho \nn\\
	& & - 2\delta \int_0^r e^{\delta(r-\rho)} \rho^{2n-2} u(\rho) d\rho
	- 2r^{2n-2} u(r) \nn\\
	&\le& 4(n-1) \int_0^r e^{\delta(r-\rho)} \rho^{2n-3} u(\rho) d\rho \nn\\
	&\le& 4(n-1) e^{\delta R} \int_0^r \rho^{2n-3} u(\rho) d\rho
	\qquad \mbox{for all } r\in (0,R),
  \eea
  because $u$ is nonnegative.\\
  Next, the Cauchy-Schwarz inequality shows that
  \be{5.8}
	2\int_0^r e^{\delta (r-\rho)} \rho^{2n-2} \sqrt{u(\rho)} g(\rho) d\rho 
	\le 2 \bigg( \int_0^R \rho^{n-1} u(\rho) d\rho \bigg)^\frac{1}{2} \cdot 
	\bigg( \int_0^r e^{2\delta(r-\rho)} \cdot \rho^{3n-3} g^2(\rho) d\rho \bigg)^\frac{1}{2}
  \ee
  for all $r\in (0,R)$, where
  \be{5.9}
	\int_0^R \rho^{n-1} u(\rho) d\rho= \frac{m}{\omega_n}
  \ee
  and
  \bea{5.10}
	\int_0^r e^{2\delta(r-\rho)} \cdot \rho^{3n-3} g^2(\rho) d\rho
	&\le& e^{2\delta R} \cdot r^{2n-2} \int_0^R \rho^{n-1} g^2(\rho) d\rho \nn\\
	&=& e^{2\delta R} \cdot r^{2n-2} \cdot \frac{\|g\|_{L^2(\Omega)}^2}{\omega_n}
	\qquad \mbox{for all } r\in (0,R)
  \eea
  with $\omega_n$ denoting the $(n-1)$-dimensional measure of $\partial B_1$.\\
  By a similar idea using pointwise estimates, the second last term in (\ref{5.6}) can be controlled according to
  \bea{5.11}
	\frac{1}{\delta} \int_0^r e^{\delta (r-\rho)} \rho^{2n-2} f^2(\rho) d\rho
	&\le& \frac{e^{\delta R}}{\delta} \cdot r^{n-1} \cdot \int_0^R \rho^{n-1} f^2(\rho) d\rho \nn\\
	&=& \frac{e^{\delta R}}{\delta} \cdot r^{n-1} \cdot \frac{\|f\|_{L^2(\Omega)}^2}{\omega_n}
	\qquad \mbox{for all } r\in (0,R).
  \eea
  Finally, upon another integration by parts we find that
  \bas
	\int_0^r e^{\delta(r-\rho)} \rho^{2n-2} (v^2)_r(\rho) d\rho
	&=& r^{2n-2} v^2(r) \nn\\
	& & - \int_0^r e^{\delta(r-\rho)} \cdot [ (2n-2) \rho^{2n-3} - \delta \rho^{2n-2}] \cdot v^2(\rho) d\rho \nn\\
	&\le& r^{2n-2} v^2(r)
	\qquad \mbox{for all } r\in (0,R),
  \eas
  because $(2n-2) \rho^{2n-3} \ge \delta \rho^{2n-2}$ for all $\rho \in (0,R)$ due to our restriction
  $\delta \le \frac{2n-2}{R}$.\\
  Combined with (\ref{5.6})-(\ref{5.11}), this shows that
  \bas
	r^{2n-2} v_r^2(r) 
	\le c_1(m) \int_0^r \rho^{2n-3} u(\rho) d\rho
	+ c_1(m) r^{n-1} \|g\|_{L^2(\Omega)}
	+ c_1(m) r^{n-1} \|f\|_{L^2(\Omega)}^2
	+ r^{2n-2} v^2(r)
  \eas
  for all $r\in (0,R)$ with some $c_1(m)>0$.
  On division by $r^{n-1}$ and integration over $r\in (0,r_0)$ we therefore obtain
  \bas
	\int_0^{r_0} r^{n-1} v_r^2(r) dr
	&\le& c_1(m) \int_0^{r_0} \frac{1}{r^{n-1}} \cdot \int_0^r \rho^{2n-3} u(\rho) d\rho dr
	+ c_1(m) r_0 \|g\|_{L^2(\Omega)} \\
	& & + c_1(m) r_0 \|f\|_{L^2(\Omega)}^2
	+ \int_0^{r_0} r^{n-1} v^2(r) dr \\
	&\le& c_1(m) \int_0^{r_0} \frac{1}{r^{n-1}} \cdot \int_0^r \rho^{2n-3} u(\rho) d\rho dr
	+c_1(m) R \|g\|_{L^2(\Omega)} \\
	& & + c_1(m) r_0 \|f\|_{L^2(\Omega)}^2
	+ \frac{c_1(m)}{\omega_n} \|v\|_{L^2(\Omega)}^2.
  \eas
  Here the Fubini theorem applies to show that
  \bea{5.12}
	\int_0^{r_0} \frac{1}{r^{n-1}} \cdot \int_0^r \rho^{2n-3} u(\rho) d\rho dr
	&=& \int_0^{r_0} \bigg( \int_\rho^{r_0} \frac{dr}{r^{n-1}} \bigg) \cdot \rho^{2n-3} u(\rho) d\rho \nn\\
	&=& \frac{1}{n-2} \int_0^{r_0} (\rho^{2-n}-r_0^{2-n}) \cdot \rho^{2n-3} u(\rho) d\rho \nn\\
	&\le& \frac{1}{n-2} \int_0^{r_0} \rho^{n-1} u(\rho) d\rho \nn\\
	&\le& \frac{m}{(n-2)\omega_n},
  \eea
  regardless of the size of $r_0 \in (0,R)$. In view of (\ref{f}) and (\ref{g}), this completes the proof.
\qed
A combination of the above two lemmata now yields an estimate for $\io |\nabla v|^2$ that is adequate for our purpose.
\begin{lem}\label{lem6}
  For all $\eps>0$ there exists $C(\eps,m,M,B,\kappa)>0$ such that each $(u,v) \in \set$ satisfies
  \be{6.1}
	\io |\nabla v|^2 
	\le \eps \io uv + C(\eps,m,M,B,\kappa) \cdot 
	\bigg( \Big\|\Delta v-v+u\Big\|_{L^2(\Omega)}^{2\theta} 
	+ \Big\|\frac{\nabla u}{\sqrt{u}}-\sqrt{u}\nabla v\Big\|_{L^2(\Omega)} +1 \bigg),
  \ee
  where
  \be{6.2}
	\theta:=\frac{1}{1+\frac{n}{(2n+4)\kappa}} \, 
	\in \Big(\frac{1}{2},1\Big).
  \ee
\end{lem}
\proof
  Let us set $\beta:=\frac{(2n+4)\kappa}{n}$, so that $\theta=\frac{\beta}{\beta+1}$. Then given $\eps \in (0,1)$,
  with notation as in (\ref{f}) and (\ref{g}) 
  we apply Lemma \ref{lem4} to $r_0:=\min\{\frac{R}{2}, \|f\|_{L^2(\Omega)}^{-\frac{2}{\beta+1}} \} \in (0,R)$
  and thus obtain $c_1=c_1(\eps,m,M,B,\kappa)>0$ such that
  \bea{6.3}
	\int_{\Omega \setminus B_{r_0}} |\nabla v|^2
	\le \frac{\eps}{4} \io uv + \frac{1}{2} \io |\nabla v|^2
	+ c_1 \cdot \Big( r_0^{-\beta} + \|f\|_{L^2(\Omega)}^\frac{2n+4}{n+4} \Big).
  \eea
  With this value of $r_0$ being fixed henceforth, Lemma \ref{lem5} provides 
  $c_2=c_2(m)>0$ 
  such that
  \bas
	\int_{B_{r_0}} |\nabla v|^2 
	\le c_2 \cdot \Big( r_0 \|f\|_{L^2(\Omega)}^2 + \|g\|_{L^2(\Omega)} + \|v\|_{L^2(\Omega)}^2 + 1 \Big),
  \eas
  whence altogether we infer that
  \bas
	\io |\nabla v|^2
	\le \frac{\eps}{2} \io uv
	+ 2 c_1 r_0^{-\beta} + 2c_1 \|f\|_{L^2(\Omega)}^\frac{2n+4}{n+4}
	+2c_2 r_0 \|f\|_{L^2(\Omega)}^2
	+2c_2 (\|g\|_{L^2(\Omega)}+1)
	+2c_2 \|v\|_{L^2(\Omega)}^2.
  \eas
  Here Lemma \ref{lem13} and (\ref{m}) say that for some $c_3=c_3(M)$ we have
  \bas
	2c_2 \|v\|_{L^2(\Omega)}^2 \le \frac{1}{2} \io |\nabla v|^2 + c_3,
  \eas
  so that
  \be{6.4}
	\io |\nabla v|^2
	\le \eps \io uv + 4c_2 (\|g\|_{L^2(\Omega)}+1) + 2c_3 + I,
  \ee
  where we abbreviate
  \bas
	I:=4c_1 r_0^{-\beta} + 4c_1 \|f\|_{L^2(\Omega)}^\frac{2n+4}{n+4} + 4c_2 r_0 \|f\|_{L^2(\Omega)}^2.
  \eas
  Now in the case $\|f\|_{L^2(\Omega)} \le (\frac{2}{R})^\frac{\beta+1}{2}$ we have $r_0=\frac{R}{2}$ and hence it follows 
  that
  \bas
	I \le 4c_1 \cdot \Big(\frac{2}{R}\Big)^\beta
	+ 4c_1 \cdot \Big(\frac{2}{R}\Big)^{\frac{\beta+1}{2} \cdot \frac{2n+4}{n+4}}
	+ 4c_2 \cdot \frac{R}{2} \cdot \Big( \frac{2}{R} \Big)^{\beta+1},
  \eas
  which clearly entails (\ref{6.1}).\\
  If, conversely, $\|f\|_{L^2(\Omega)} > (\frac{2}{R})^\frac{\beta+1}{2}$ and thus 
  $r_0=\|f\|_{L^2(\Omega)}^{-\frac{2}{\beta+1}}$ then 
  \bas
	I &\le&
	4 c_1 \|f\|_{L^2(\Omega)}^\frac{2\beta}{\beta+1}
	+4c_1 \|f\|_{L^2(\Omega)}^\frac{2n+4}{n+4}
	+ 4c_2 \|f\|_{L^2(\Omega)}^{2-\frac{2}{\beta+1}} \\
	&=& 4(c_1+c_2) \|f\|_{L^2(\Omega)}^\frac{2\beta}{\beta+1} 
	+4c_1 \|f\|_{L^2(\Omega)}^\frac{2n+4}{n+4}.
  \eas
  Since $\kappa>n-2$ implies
  \bas
	\frac{\beta}{\frac{n+2}{2}}
	=\frac{2}{n+2} \cdot \frac{(2n+4)\kappa}{n}
	>\frac{4(n-2)}{n} \ge \frac{4}{3}>1
  \eas
  due to the fact that $n\ge 3$, it can easily be checked that $2\theta\equiv \frac{2\beta}{\beta+1}>\frac{2n+4}{n+4}$.
  This verifies the inclusion in (\ref{6.2}), and furthermore Young's inequality yields $c_4=c_4(\eps, m, M, B, \kappa)>0$ 
  such that
  \bas
	I \le 4(c_1+c_2+1) \|f\|_{L^2(\Omega)}^\frac{2\beta}{\beta+1} + c_4,
  \eas
  so that (\ref{6.4}) shows that (\ref{6.1}) is also valid when 
  $\|f\|_{L^2(\Omega)} > (\frac{2}{R})^\frac{\beta+1}{2}$.
\qed
Thereupon, the main result of this section actually reduces to a corollary.\abs
\proofc of Lemma \ref{lem8}. \quad
  We only need to start from Lemma \ref{lem7} and then apply Lemma \ref{lem6} to $\eps:=\frac{1}{4}$ to obtain
  $c_1=c_1(M)>0$ and $c_2=c_2(m,M,B,\kappa)>0$ such that with $f$ and $g$ as in (\ref{f}) and (\ref{g}) we have
  \bas
	\io uv &\le& \frac{1}{2} \io uv + c_1 \cdot \Big( \|f\|_{L^2(\Omega)}^\frac{2n+4}{n+4} +1 \Big) 
	+ c_2 \cdot \Big( \|f\|_{L^2(\Omega)}^{2\theta} + \|g\|_{L^2(\Omega)} + 1 \Big).
  \eas
  Since $\frac{2n+4}{n+4}<2\theta$ by (\ref{6.2}), using Young's inequality we immediately arrive at (\ref{8.1}).
\qed
\mysection{Blow-up. Proof of Theorem \ref{theo12}}\label{sect_blowup}
We now plan to apply the above estimates to solutions of the dynamical problem (\ref{0}).
Still referring to the definition of $\set$ as introduced in the beginning of Section \ref{sect_estimates},
we first turn the outcome of Lemma \ref{lem8} into a statement of type (\ref{bound}), still for arbitrary
functions in $\set$.
\begin{theo}\label{theo9}
  There exists $C(m,M,B,\kappa)>0$ such that for all $(u,v) \in \set$ we have
  \be{9.1}
	\F(u,v) \ge -C(m,M,B,\kappa) \cdot \Big( \D^\theta (u,v) +1 \Big)
  \ee
  with $\theta \in (\frac{1}{2},1)$ given by (\ref{6.2}).
\end{theo}
\proof
  Since $\theta>\frac{1}{2}$, we may apply Young's inequality to (\ref{8.1}) to find $c_1=c_1(m,M,B,\kappa)>0$ such that
  \bas
	\io uv \le c_1 \bigg( \Big(\|f\|_{L^2(\Omega)}^2 + \|g\|_{L^2(\Omega)}^2 \Big)^\theta +1 \bigg)
  \eas
  with $f$ and $g$ as given by (\ref{f}) and (\ref{g}). 
  Therefore, using that $\xi\ln \xi \ge -\frac{1}{e}$ for all $\xi>0$ we obtain the inequality
  \bas
	\F(u,v)
	&=& \frac{1}{2} \io |\nabla v|^2 + \frac{1}{2} \io v^2 - \io uv + \io u\ln u \\
	&\ge& - \io uv - \frac{|\Omega|}{e} \\
	&\ge& - c_2 \cdot \Big( (\|f\|_{L^2(\Omega)}^2 + \|g\|_{L^2(\Omega)}^2)^\theta +1 \Big)
  \eas
  with $c_2 \equiv c_2(m,M,B,\kappa):=c_1+\frac{|\Omega|}{e}$.
  Since by definition of $f$ and $g$ we have $\D(u,v)=\|f\|_{L^2(\Omega)}^2 + \|g\|_{L^2(\Omega)}^2$,
  this already establishes (\ref{9.1}).
\qed
Now given a solution $(u,v)$ of (\ref{0}), 
the fact that in (\ref{9.1}) we have $\theta<1$ will enable us to derive
an ODI for $t\mapsto - \F(u(\cdot,t),v(\cdot,t))$ with superlinearly growing nonlinearity. 
For initial data $(u_0,v_0)$ with large negative energy $\F(u_0,v_0)$, this means that 
$(u,v)$ cannot exist globally.
\begin{lem}\label{lem11}
  Let $m>0, A>0$ and $\kappa>n-2$. Then there exist $K=K(m,A,\kappa)>0$ and $C=C(m,A,\kappa)>0$ 
  such that for each $(u_0,v_0)$ from the set
  \bea{tB}
	\widetilde{\B}(m,A,\kappa)
	&:=& \bigg\{
	(u_0,v_0) \in C^0(\bar\Omega) \times W^{1,\infty}(\Omega) \ \bigg| \
	\mbox{$u_0$ and $v_0$ are radially symmetric and positive in $\bar\Omega$} \nn\\
	& & \hspace*{35mm}
	\mbox{with $\io u_0=m$, $\|v_0\|_{W^{1,2}(\Omega)} \le A$ and $\F(u_0,v_0) \le -K$} \bigg\},
  \eea
  the corresponding solution $(u,v)$ of (\ref{0}) has the property
  \be{11.11}
	\F(u(\cdot,t),v(\cdot,t))
	\le \frac{\F (u_0,v_0)}{(1-Ct)^\frac{\theta}{1-\theta}}
	\qquad \mbox{for all } t\in (0,\tm),
  \ee
  where $\theta \in (\frac{1}{2},1)$ is as given by (\ref{6.2}).\\
  In particular, for any such solution we have $\tm<\infty$, that is, $(u,v)$ blows up in finite time.
\end{lem}
\proof
  Let us fix $c_1>0$ such that
  \be{11.111}
	\|\varphi\|_{L^1(\Omega)} \le c_1 \|\varphi\|_{W^{1,2}(\Omega)}
	\qquad \mbox{for all } \varphi \in W^{1,2}(\Omega).
  \ee
  According to Corollary \ref{cor3}, we can pick $c_2=c_2(\kappa)>0$ such that whenever $u_0\in C^0(\bar\Omega)$
  and $v_0 \in W^{1,\infty}(\Omega)$ are radial and positive, the corresponding solution $(u,v)$ of (\ref{0}) satisfies
  \be{11.1}
	v(r,t) \le c_2 \cdot \Big( \|u_0\|_{L^1(\Omega)} + \|v_0\|_{L^1(\Omega)} + \|\nabla v_0\|_{L^2(\Omega)} \Big)
	\cdot r^{-\kappa} \qquad \mbox{for all } (r,t) \in (0,R) \times (0,\tm).
  \ee
  Next, writing $B:=c_2 (m+c_1 A+A)$ and $M:=\max \{m,c_1 A\}$ we invoke Theorem \ref{theo9} to obtain
  $c_3=c_3(m,M,B,\kappa)>0$ such that
  \be{11.2}
	\F (\tilde u, \tilde v) \ge - c_3 \cdot \Big(\D^\theta(\tilde u,\tilde v) +1 \Big)
	\qquad \mbox{for all } (\tilde u,\tilde v) \in \set.
  \ee
  We will see that then (\ref{11.11}) holds for all $(u_0,v_0) \in \widetilde{B}(m,A,\kappa)$ if we define
  \be{11.3}
	K(m,A,\kappa) := 2c_3
  \ee
  and
  \be{11.4}
	C(m,A,\kappa):= \frac{1-\theta}{2c_3 \theta}.
  \ee
  Indeed, given $(u_0,v_0) \in \widetilde{B}(m,A,\kappa)$ we know from (\ref{11.1}) and (\ref{11.111}) that
  the solution $(u,v)$ of (\ref{0}) emanating from $(u_0,v_0)$ is smooth and radially symmetric with 
  $u>0$ and $v>0$ in $\bar\Omega \times [0,\tm)$ and
  \bas
	v(r,t) \le c_2 \cdot (m+c_1 A + A) \cdot r^{-\kappa} = B r^{-\kappa}
	\qquad \mbox{for all } (r,t) \in (0,R) \times (0,\tm).
  \eas
  Since moreover $\io u(\cdot,t) \equiv \io u_0=m$ and $v(\cdot,t) \le \max \{\io u_0,\io v_0\}
  \le \max \{m, c_1 A\}=M$ for all $t\in (0,\tm)$ by (\ref{mass}) and (\ref{mass_v}), it follows that $(u(\cdot,t),v(\cdot,t)) \in \set$
  for all $t\in (0,\tm)$ and hence (\ref{11.2}) may be applied to $(\tilde u,\tilde v):=(u(\cdot,t),v(\cdot,t))$
  for any such $t$.\\
  In order to derive (\ref{11.11}) from this and the energy inequality (\ref{energy}), let us make sure that
  \bas
	y(t):=-\F(u(\cdot,t),v(\cdot,t)), \qquad t\in [0,\tm),
  \eas
  defines a positive function $y\in C^0([0,\tm)) \cap C^1((0,\tm))$ which satisfies
  \be{11.5}
	y'(t) \ge c_4 y^\frac{1}{\theta}(t)
	\qquad \mbox{for all } t\in (0,\tm)
  \ee
  with $c_4=c_4(m,M,B,\kappa):=(2c_3)^{-\frac{1}{\theta}}$.\\
  In fact, the claimed regularity properties of $y$ immediately result from those of $(u,v)$, whereas (\ref{energy})
  ensures that $y$ is nondecreasing and thus
  \be{11.6}
	y(t) \ge y(0) \ge K(m,A,\kappa) =2c_3>0
	\qquad \mbox{for all } t\in (0,\tm).
  \ee
  Therefore we may invert so as to obtain from (\ref{11.2}) and (\ref{11.6}) that
  \bas
	\D^\theta(u(\cdot,t),v(\cdot,t)) 
	\ge \frac{y(t)}{c_3} -1 
	\ge \frac{y(t)}{c_3} - \frac{y(t)}{2c_3} = \frac{y(t)}{2c_3}
	\qquad \mbox{for all } t\in (0,\tm).
  \eas
  In light of (\ref{energy}), $y$ thus satisfies
  \bas
	y'(t) \ge \D(u(\cdot,t),v(\cdot,t)) \ge \Big( \frac{y(t)}{2c_3} \Big)^\frac{1}{\theta}
	\qquad \mbox{for all } t\in (0,\tm),
  \eas
  which precisely yields (\ref{11.5}).\\
  Now by straightforward integration, we see that
  \bas
	y(t) \ge y(0) \cdot \Big\{ 1 - \frac{1-\theta}{\theta} c_4 y^\frac{1-\theta}{\theta}(0) \cdot t 
	\Big\}^{-\frac{\theta}{1-\theta}}
	\qquad \mbox{for all } t\in (0,\tm),
  \eas
  which implies (\ref{11.11}) upon the observation that (\ref{11.6}) entails
  \bas
	\frac{1-\theta}{\theta} c_4 \cdot y^\frac{1-\theta}{\theta}(0)
	\ge \frac{1-\theta}{\theta} c_4 \cdot (2c_3)^\frac{1-\theta}{\theta}
	= \frac{1-\theta}{\theta} \cdot \frac{1}{2c_3}
	= C(m,A,\kappa).
  \eas
  The proof is complete.
\qed
In a last step we can remove the auxiliary parameter $\kappa$ appearing in the latter result so as to obtain
the blow-up criterion stated in Theorem \ref{theo12}.\abs
\proofc of Theorem \ref{theo12}. \quad
  We only need to fix an arbitrary $\kappa>n-2$ and apply Lemma \ref{lem11} to find that the conclusion holds
  if we let $K(m,A):=K(m,A,\kappa)$ and $T(m,A):=\frac{1}{C(m,A,\kappa)}$ with $K(m,A,\kappa)$ and
  $C(m,A,\kappa)$ as given by Lemma \ref{lem11}.
\qed
\mysection{A density property of $\B (m,A)$. Proof of Theorem \ref{theo15}}\label{sect_dense}
The proof of Theorem \ref{theo15} will be an immediate consequence of the following lemma which roughly says that 
arbitrary initial data can be approximated by low-energy initial data in the claimed manner.
\begin{lem}\label{lem14}
  Let $m>0$ and $u \in C^0(\bar\Omega)$ and $v\in W^{1,\infty}(\Omega)$ be radially symmetric and positive in $\bar{\Omega}$
  with $\io u=m$. Then for each $p\in (1,\frac{2n}{n+2})$ there exist sequences $(u_k)_{k\in\N} \subset C^0(\bar\Omega)$
  and $(v_k)_{k\in\N} \subset W^{1,\infty}(\Omega)$ of radially symmetric positive functions satisfying $\io u_k=m$ for all
  $k\in\N$ and 
  \be{14.1}
	u_k\to u \quad \mbox{in } L^p(\Omega)
	\quad \mbox{and} \quad
	v_k \to v \quad \mbox{in } W^{1,2}(\Omega)
	\qquad \mbox{as } k\to\infty,
  \ee
  but such that with $\F$ as defined in (\ref{F}) we have
  \be{14.2}
	\F(u_k,v_k) \to -\infty
	\qquad \mbox{as } k\to\infty.
  \ee
\end{lem}
\proof
  We fix an arbitrary sequence $(r_k)_{k\in\N} \subset (0,R)$ such that $r_k\to 0$ as $k\to\infty$, and let
  \bas
	\varphi(\xi) := \int_0^1 \rho^{n-1} (\rho^2+\xi)^{-\frac{n}{2}} d\rho, \qquad \xi>0.
  \eas
  Then by monotone convergence we have $\varphi(\xi) \nearrow \infty$ as $\xi\searrow 0$, so that for each $k\in\N$
  it is possible to fix $\eta_k \in (0,R^2)$ appropriately small such that
  \be{14.3}
	r_k^n \cdot \varphi \Big( \frac{\eta_k}{r_k^2} \Big) \ge k.
  \ee
  Now given $p\in (1,\frac{2n}{n+2})$, we can fix $\alpha>0$ fulfilling $n-\frac{n}{p} < \alpha < \frac{n-2}{2}$
  and thereupon introduce positive radial functions $\tilde u_k, u_k$ and $v_k$ on $\bar\Omega$ by defining
  \bas
	\tilde u_k(r):= \left\{ \begin{array}{ll}
	a_k \cdot (r^2+\eta_k)^{-\frac{n-\alpha}{2}}, \qquad & r \in [0,r_k], \\[1mm]
	u(r), & r \in (r_k,R],
	\end{array} \right.
  \eas
  and
  \be{14.4}
	u_k:=\frac{m\tilde u_k}{\|\tilde u_k\|_{L^1(\Omega)}}
  \ee
  as well as
  \bas
	v_k(r):= \left\{ \begin{array}{ll}
	b_k \cdot (r^2+\eta_k)^{-\frac{\alpha}{2}}, \qquad & r \in [0,r_k], \\[1mm]
	v(r), & r \in (r_k,R],
	\end{array} \right.
  \eas
  with
  \bas
	a_k:=(r_k^2+\eta_k)^\frac{n-\alpha}{2} \cdot u(r_k)
	\qquad \mbox{and} \qquad
	b_k:=(r_k^2+\eta_k)^\frac{\alpha}{2} \cdot v(r_k)
  \eas
  for $k\in\N$.
  Then clearly $v_k \in W^{1,\infty}(\Omega)$, whereas $\tilde u_k$ and $u_k$ belong to $C^0(\bar{\Omega})$ with $\io u_k=m$.
  Moreover, once again writing $\omega_n:=|\partial B_1|$ we have 
  \bea{14.5}
	\| \tilde u_k \|_{L^p(B_{r_k})}^p
	&=& \omega_n \cdot \int_0^{r_k} r^{n-1} \cdot a_k^p (r^2+\eta_k)^{-\frac{(n-\alpha)p}{2}} dr \nn\\
	&=& \omega_n \cdot u^p(r_k) \cdot 
		\int_0^{r_k} r^{n-1} \cdot \Big(\frac{r_k^2+\eta_k}{r^2+\eta_k}\Big)^\frac{(n-\alpha)p}{2} dr \nn\\
	&\le& \omega_n \cdot u^p(r_k) \cdot \int_0^{r_k} r^{n-1-(n-\alpha)p} dr \nn\\
	&=& \frac{\omega_n}{n-(n-\alpha)p} \cdot u^p(r_k) r_k^{n-(n-\alpha)p}
	\qquad \mbox{for all } k\in\N,
  \eea
  because $\alpha>n-\frac{n}{p}$ implies $(n-\alpha)p<n$. 
  Since $r_k\to 0$ and $u$ is bounded, we thus infer from 
  \be{14.51}
    \|\tilde u_k-u\|_{L^p(\Omega)} = \|\tilde u_k-u\|_{L^p(B_{r_k})} \le \| \tilde u_k \|_{L^p(B_{r_k})} 
    + \| u \|_{L^\infty (\Omega)} |B_{r_k}|^{\frac{1}{p}}
  \ee 
  that $\tilde u_k \to u$ in $L^p(\Omega)$ as $k\to\infty$.
  In particular, this entails that $\tilde u_k \to u$ also in $L^1(\Omega)$ and accordingly 
  \be{14.55}
	\|\tilde u_k\|_{L^1(\Omega)} \to m
	\qquad \mbox{as } k\to\infty,
  \ee
  so that combining (\ref{14.4}) and (\ref{14.51}) yields
  \bas
	u_k \to u \qquad \mbox{in } L^p(\Omega)
	\qquad \mbox{as } k\to\infty.
  \eas
  Next, computing $v_{kr}=-\alpha b_k r(r^2+\eta_k)^{-\frac{\alpha+2}{2}}$ for $r\in (0,r_k)$ and observing
  that $b_k \le c_1:=(2R^2)^\frac{\alpha}{2} \|v\|_{L^\infty(\Omega)}$ for all $k\in\N$, we find that
  \bas
	\| \nabla v_k \|_{L^2(B_{r_k})}^2
	&=&  \omega_n \cdot \alpha^2 b_k^2 \cdot \int_0^{r_k} r^{n+1} (r^2+\eta_k)^{-\alpha-2} dr \\
	&\le& \omega_n \cdot \alpha^2 c_1^2 \cdot \int_0^{r_k} r^{n-2\alpha-3} dr \\
	&=& \frac{\omega_n \alpha^2 c_1^2}{n-2\alpha-2} \cdot r_k^{n-2\alpha-2}
	\qquad \mbox{for all } k\in\N
  \eas
  and, similarly,
  \bas
	\| v_k \|_{L^2(B_{r_k})}^2
	&=& \omega_n \cdot b_k^2 \cdot \int_0^{r_k} r^{n-1} (r^2+\eta_k)^{-\alpha} dr
	\le \frac{\omega_n c_1^2}{n-2\alpha} \cdot r_k^{n-2\alpha}
	\qquad \mbox{for all } k\in\N,
  \eas
  for our restriction $\alpha<\frac{n-2}{2}$ ensures that both $n-2\alpha-2$ and $n-2\alpha$ are positive. 
  Again since $r_k\to 0$, in view of the estimate 
  \bas 
    \|v_k-v\|_{W^{1,2}(\Omega)} = \|v_k-v\|_{W^{1,2}(B_{r_k})} \le \|v_k\|_{W^{1,2}(B_{r_k})}
    + \|v\|_{W^{1,\infty}(\Omega)} (2 |B_{r_k}|)^{\frac{1}{2}}
  \eas 
  and the boundedness of $v$ in $W^{1,\infty}(\Omega)$ this shows that
  that
  \bas
	v_k \to v \quad \mbox{in } W^{1,2}(\Omega)
	\qquad \mbox{as } k\to\infty
  \eas
  and thus completes the proof of (\ref{14.1}).\\
  To verify (\ref{14.2}), we first note that since $\sup_{\xi>0} \frac{\xi \ln \xi}{\xi^p}$ is finite thanks to our assumption
  $p>1$, from the boundedness of $(u_k)_{k\in\N}$ in $L^p(\Omega)$ and of $(v_k)_{k\in\N}$ in $W^{1,2}(\Omega)$
  asserted by (\ref{14.1}) we obtain $c_2>0$ such that
  \be{14.6}
	\frac{1}{2} \io |\nabla v_k|^2 + \frac{1}{2} \io v_k^2 + \io u_k \ln u_k \le c_2
	\qquad \mbox{for all } k\in\N.
  \ee
  On the other hand, using that 
  \bas
	a_k b_k = u(r_k) v(r_k) \cdot (r_k^2+\eta_k)^\frac{n}{2} \ge u(r_k)v(r_k) \cdot r_k^n
	\qquad \mbox{for all } k\in\N
  \eas
  and recalling the definition of $\varphi$ we have
  \bas
	\io u_k v_k &\ge& \int_{B_{r_k}} u_k v_k \\
	&=& \omega_n \cdot \frac{m}{\|\tilde u_k\|_{L^1(\Omega)}} \cdot a_k b_k \cdot 
		\int_0^{r_k} r^{n-1} (r^2+\eta_k)^{-\frac{n}{2}} dr \\
	&=& \omega_n \cdot \frac{m}{\|\tilde u_k\|_{L^1(\Omega)}} \cdot a_k b_k \cdot 
		\int_0^1 \rho^{n-1} \cdot \Big(\rho^2 + \frac{\eta_k}{r_k^2}\Big)^{-\frac{n}{2}} d\rho \\
	&\ge& \omega_n \cdot \frac{m}{\|\tilde u_k\|_{L^1(\Omega)}} \cdot u(r_k) v(r_k) r_k^n 
		\cdot \varphi \Big( \frac{\eta_k}{r_k^2}\Big)
	\qquad \mbox{for all } k\in\N.
  \eas
  Hence, (\ref{14.3}) warrants that
  \bas
	\io u_k v_k \ge \omega_n \cdot \frac{m}{\|\tilde u_k\|_{L^1(\Omega)}} \cdot u(r_k)v(r_k) \cdot k
	\qquad \mbox{for all } k\in\N,
  \eas
  so that from (\ref{14.55}) we infer that
  \bas
	\liminf_{k\to\infty} \frac{1}{k} \cdot \io u_k v_k \ge \omega_n \cdot u(0)v(0)>0
	\qquad \mbox{as } k\to\infty,
  \eas
  because $u$ and $v$ are continuous. Thus, by positivity of $u$ and $v$ we conclude that $\io u_k v_k \to\infty$ as 
  $k\to\infty$, and in conjunction with (\ref{14.6}) this proves (\ref{14.2}).
\qed
\proofc of Theorem \ref{theo15}.\quad
In view of Theorem \ref{theo12}, the claim directly results from Lemma \ref{lem14}.
\qed
{\bf Acknowledgment.} \quad
The author would like to thank Christian Stinner for numerous useful remarks which significantly improved
this work.
\end{document}